\documentclass[psamsfonts,reqno]{amsart}
\usepackage{amssymb,eucal,graphics,xy}
\usepackage{epstopdf}
\xyoption{all}
\xyoption{ps}

\numberwithin{equation}{section}

\theoremstyle{plain}

\newtheorem{lemma}{Lemma}[section]
\newtheorem{theorem}[lemma]{Theorem}
\newtheorem{proposition}[lemma]{Proposition}
\newtheorem{corollary}[lemma]{Corollary}

\newtheorem*{stat}{\name}
\newcommand{\name}{testing}

\theoremstyle{definition}
\newtheorem{definition}[lemma]{Definition}

\newtheorem{problem}{Problem}

\theoremstyle{remark}

\newenvironment{all}[1]{\renewcommand{\name}{#1}\begin{stat}}
                        {\end{stat}}

\newcommand{\qedc}{{\qed}~{\rm Claim~{\theclaim}.}}

\newcommand{\cls}[2]{{[{#1}]_{#2}}}
\newcommand{\res}{\mathbin{\restriction}}

\newcommand{\distr}{dis\-trib\-u\-tive}
\newcommand{\distry}{dis\-trib\-u\-tiv\-i\-ty}
\newcommand{\eps}{\varepsilon}

\newcommand{\set}[1]{\{{#1}\}}
\newcommand{\setm}[2]{\set{{#1}\mid{#2}}}
\newcommand{\Set}[1]{\left\{{#1}\right\}}
\newcommand{\Setm}[2]{\Set{{#1}\mid{#2}}}
\newcommand{\seq}[1]{\langle{#1}\rangle}
\newcommand{\seqm}[2]{\seq{{#1}\mid{#2}}}
\newcommand{\Seq}[1]{\left\langle{#1}\right\rangle}
\newcommand{\Seqm}[2]{\Seq{{#1}\mid{#2}}}
\newcommand{\famm}[2]{\left({#1}\mid{#2}\right)}
\newcommand{\into}{\hookrightarrow}
\newcommand{\onto}{\twoheadrightarrow}
\newcommand{\ol}[1]{\overline{#1}}
\newcommand{\conj}{\mathbin{\bigwedge\mkern-15mu\bigwedge}}

\newcommand{\CC}{\boldsymbol{\mathcal{C}}}
\newcommand{\NC}{\boldsymbol{\mathcal{NC}}}

\newcommand{\vx}{\mathsf{x}}
\newcommand{\vy}{\mathsf{y}}

\newcommand{\bv}[1]{\Vert{#1}\Vert}

\DeclareMathOperator{\Hom}{Hom}
\DeclareMathOperator{\End}{End}
\DeclareMathOperator{\Sub}{Sub}
\DeclareMathOperator{\Subc}{Sub_c}
\DeclareMathOperator{\Id}{Id}

\DeclareMathOperator{\Dim}{Dim}
\DeclareMathOperator{\cen}{cen}

\DeclareMathOperator{\Clop}{Clop}

\DeclareMathOperator{\At}{At}

\newcommand{\intp}[1]{\lfloor{#1}\rfloor}
\newcommand{\Cont}{\mathbf{C}}
\newcommand{\BS}{\mathbf{S}}
\newcommand{\BL}{\mathbf{L}}
\newcommand{\CL}{\mathbf{CL}}

\newcommand{\DD}{\Delta}

\newcommand{\two}{\mathbf{2}}

\DeclareMathOperator{\Mat}{M}
\DeclareMathOperator{\GL}{GL}

\newcommand{\bI}{{\boldsymbol{I}}}

\newcommand{\bL}{{\boldsymbol{L}}}

\newcommand{\cF}{\mathcal{F}}

\newcommand{\cU}{\mathcal{U}}

\newcommand{\cL}{\mathcal{L}}

\DeclareMathOperator{\im}{im}

\newcommand{\Pow}{\mathfrak{P}}
\newcommand{\PP}{\mathbb{P}}
\newcommand{\FF}{\mathbb{F}}
\newcommand{\NN}{\mathbb{N}}

\newcommand{\les}{\leqslant}
\newcommand{\ges}{\geqslant}

\newcommand{\fp}{\mathfrak{p}}

\author[F. Wehrung]{Friedrich Wehrung}
\address{LMNO, CNRS UMR 6139\\
         Universit\'e de Caen, Campus II\\
         D\'epartement de Math\'ematiques\\
         B.P. 5186\\
         14032 CAEN Cedex\\
         FRANCE}
 \email{wehrung@math.unicaen.fr}
 \urladdr{http://www.math.unicaen.fr/\~{}wehrung}

\date{\today}

\begin{document}

\title[Von Neumann coordinatization]
{Von Neumann coordinatization is not first-order}

\begin{abstract}
A lattice $L$ is \emph{coordinatizable}, if it
is isomorphic to the lattice~$\bL(R)$ of principal right ideals of
some von~Neumann regular ring~$R$. This forces~$L$ to be
complemented modular. All known sufficient conditions
for coordinatizability, due first to J. von~Neumann, then to
B.~J\'onsson, are first-order. Nevertheless, we prove that
coordinatizability of lattices is not
first-order, by finding a non-coordinatizable lattice $K$ with a
coordinatizable countable elementary extension $L$. This solves a
1960 problem of B.~J\'onsson. We also prove that there is no
$\cL_{\infty,\infty}$ statement equivalent to coordinatizability.
Furthermore, the class of coordinatizable lattices is not closed
under countable directed unions; this solves another problem of
B.~J\'onsson from 1962.
\end{abstract}

\maketitle

\section{Introduction}\label{S:Intro}
A \emph{coordinatization theorem} is a statement that expresses a
class of geometric objects in algebraic terms. Hence it is a path
from \emph{synthetic} geometry to \emph{analytic} geometry. While
the former includes lattice theory, as, for example, abundantly
illustrated in M.\,K. Bennett's survey paper \cite{Benn95}, the
latter is more often written in the language of rings and modules.
Nevertheless the concepts of analytic and synthetic geometry will
not let themselves be captured so easily. For example, the main
result of E. Hrushovski and B. Zilber \cite[Theorem~A]{HrZi} may
certainly be viewed as a coordinatization theorem, with geometric
objects of \emph{topological} nature.

It should be no surprise that coordinatization theorems are usually
very difficult results. The classical coordinatization theorem of
Arguesian affine planes (as, for example, presented in E. Artin
\cite[Chapitre~II]{Artin}) was extended over the last century to
a huge work on modular lattices, which also brought surprising and
deep connections with coordinatization results in universal algebra,
see the survey paper by C. Herrmann \cite{Herr95}.
We cite the following milestone, due to J. von~Neumann~\cite{Neum60}.

\begin{all}{Von~Neumann's Coordinatization Theorem}
If a complemented modular lattice $L$ has a spanning
finite homogeneous sequence with at least four elements, then~$L$ is
\emph{coordinatizable}, that is, there exists a von~Neumann regular
ring $R$ such that $L$ is isomorphic to the lattice $\bL(R)$ of all
principal right ideals of $R$.
\end{all}

We refer the reader to
Sections~\ref{S:BasicLatt}--\ref{S:BasicModules} for precise
definitions. We observe that while the statement that a lattice is
coordinatizable is, apparently, ``complicated'' (it begins with
an existential quantifier over regular rings), von~Neumann's
sufficient condition is logically simple---for example, having a
spanning homogeneous sequence with four elements is a first-order
condition.

The strongest known extension of von~Neumann's Coordinatization
Theorem is due to B.~J\'onsson \cite{Jons60}. For further
coordinatization results of modular lattices, see, for
example, B. J\'onsson and G. Monk \cite{JoMo69}, A. Day and D.
Pickering \cite{DaPi83}, or the survey by M. Greferath and S.\,E.
Schmidt \cite{GrSc}.

\begin{all}{J\'onsson's Extended Coordinatization Theorem}
\label{T:JonssThm}
Every complemented Arguesian lattice $L$ with a ``large partial
$3$-frame'' is coordinatizable.
\end{all}

Although having a large partial $3$-frame is, apparently, not a
first-order condition, we prove in Section~\ref{S:LP3Frme}, by using
the dimension monoid introduced in F. Wehrung \cite{WDim}, that it
can be expressed by a single first-order sentence.

Is coordinatizability first-order? The question was raised
by B. J\'onsson in the Introduction of \cite{Jons60}. We quote the
corresponding excerpt.

\begin{quote}\em
A complete solution to our problem would consist in an axiomatic
characterization of the class of all coordinatizable
lattices. However, this seems to be an extremely difficult problem,
and in fact it is doubtful that any reasonable axiom system can be
found.
\end{quote}

In the present paper, we confirm J\'onsson's negative
guess, in particular removing the word ``reasonable'' from the
second sentence above. In fact, our negative solution even holds for
a restricted class of complemented modular lattices, namely, those
that satisfy the identity of \emph{$2$-\distry},
 \[
 x\vee(y_0\wedge y_1\wedge y_2)=
 \bigl(x\vee(y_0\wedge y_1)\bigr)\wedge
 \bigl(x\vee(y_0\wedge y_2)\bigr)
 \vee\bigl(x\vee(y_1\wedge y_2)\bigr).
 \]
A few examples of $2$-distributive modular lattices are diagrammed on
Figure~1. The subspace lattice of a three-dimensional vector space
is not $2$-\distr. An important characterization of $2$-\distry\ for
modular lattices is provided by C. Herrmann, D. Pickering, and M.
Roddy \cite{HPR94}: \emph{A modular lattice is $2$-\distr\ if{f} it
can be embedded into the subspace lattice of a vector space over any
field}. So, in some sense, the theory of $2$-\distr\ modular lattices
is the ``characteristic-free'' part of the theory of modular
lattices.

\begin{figure}[hbt]\label{Fig:2distr}
\includegraphics{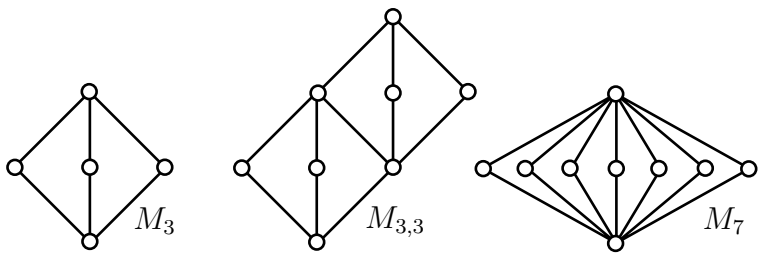}
\caption{A few $2$-distributive modular lattices.}
\end{figure}

We also find a large proper class of non-coordinatizable lattices
with spanning~$M_{\omega}$, see Theorem~\ref{T:LnotCoord}. This
result is sufficient to ensure that there is no
$\cL_{\infty,\infty}$ statement equivalent to coordinatizability (see
Theorem~\ref{T:NotLinf}).

We put $\NN=\omega\setminus\set{0}$, and we denote
by $\PP$ the set of all prime numbers. For a prime $p$, we denote by
$\FF_{p^{\infty}}$ an (the) algebraic closure of the prime field
$\FF_p$ with $p$ elements, and we put
$\FF_q=\setm{x\in\FF_{p^{\infty}}}{x^q=x}$, for any power $q$ of
$p$. Hence $\FF_q$ is a (the) field with $q$ elements, and $\FF_m$
is a subfield of $\FF_n$ if{f} $n$ is a power of $m$.

Following standard set-theoretical notation,
we denote by $\omega$ the chain of all natural numbers and by
$\omega_1$ the first uncountable ordinal.

If $\alpha$ is an equivalence relation on a set $X$, we denote by
$\cls{x}{\alpha}$ the $\alpha$-equivalence class of $x$ modulo
$\alpha$, for every $x\in X$. If $f\colon X\to Y$ is a map, we put
$f[Z]=\setm{f(x)}{x\in Z}$, for any $Z\subseteq X$. For an infinite
set $I$, a family $x=\seqm{x_i}{i\in I}$ is \emph{almost constant},
if there exists a (necessarily unique) $a$ such that
$\setm{i\in I}{x_i\neq a}$ is finite, and then we put $a=x(\infty)$,
the \emph{limit} of $x$.

\section{Lattices}\label{S:BasicLatt}
Standard textbooks on lattice theory are G. Birkhoff \cite{Birk79},
G. Gr\"at\-zer~\cite{Grat98}, and R.\,N. McKenzie, G.\,F. McNulty,
and W.\,F. Taylor \cite{MMTa87}. We say that a lattice~$L$ is
\emph{bounded}, if it has a zero (i.e., a least element), generally
denoted by $0$, and a unit (i.e., a largest element), generally
denoted by $1$. For an element $a$ in a lattice~$L$,  we put
 \[
 L\res a=\setm{x\in L}{x\leq a},\text{ the \emph{principal ideal}
 generated by }a.
 \]
We say that $L$ is \emph{modular}, if it satisfies
the identity
 \[
 x\wedge(y\vee(x\wedge z))=(x\wedge y)\vee(x\wedge z).
 \]
We shall sometimes mention a stronger identity than modularity,
the so-called \emph{Arguesian identity}, which can be found, for
example, in \cite[Section~IV.4]{Grat98}. The Arguesian identity holds
in every lattice of permuting equivalence relations (see
B. J\'onsson \cite{Jons53} or \cite[Section~IV.4]{Grat98}). In
particular, it holds in the lattice $\Sub M$ of all submodules of
any right module $M$ over any ring.

In case $L$ has a zero and for $a$, $b$, $c\in L$, we let
$c=a\oplus b$ hold, if $c=a\vee b$ and $a\wedge b=0$. It is
well-known that the partial operation $\oplus$ is associative in
case~$L$ is modular (see \cite[Satz~I.1.8]{FMae58} or
\cite[Proposition~8.1]{WDim}).

We say that $L$ is \emph{complemented}, if it is bounded and every
$x\in L$ has a complement, that is, an element $y\in L$ such that
$x\oplus y=1$. We say that $L$ is \emph{sectionally complemented}, if
$L$ has a zero and every principal ideal of $L$ is a complemented
sublattice.

In case $L$ has a zero,
the relations of \emph{perspectivity}, $\sim$, and
\emph{subperspectivity}, $\lesssim$, are defined in $L$ by
 \begin{align*}
 x\sim y,&\text{ if }\exists z\in L\text{ such that }
 x\oplus z=y\oplus z,\\
 x\lesssim y,&\text{ if }\exists z\in L\text{ such that }
 x\oplus z\leq y\oplus z.
 \end{align*}
In case $L$ is sectionally complemented and modular, $x\lesssim y$
if{f} there exists $x'\leq y$ (resp., $y'\geq x$) such that
$x\sim x'$ (resp., $y\sim y'$) (see \cite[Theorem I.3.1]{Neum60}).

An element $a$ in a lattice $L$ is \emph{neutral} (see
\cite[Section~III.3]{Grat98}), if $\set{a,x,y}$
generates a distributive sublattice of $L$, for all $x$, $y\in L$.
An ideal $I$ of a lattice $L$ is \emph{neutral} (see
\cite[Section~III.4]{Grat98}), if $I$ is a neutral element of the
lattice $\Id L$ of ideals of $L$. Then $I$ is a
\emph{distributive ideal} of $L$, that is, the equivalence relation
$\equiv_I$ of $L$ defined by
 \[
 x\equiv_Iy\Longleftrightarrow\exists u\in I\text{ such that }
 x\vee u=y\vee u,\text{ for all }x,\,y\in L,
 \]
is a congruence of $L$. Then we denote by $L/I$ the quotient lattice
$L/{\equiv_I}$, and we put $\cls{x}{I}=\cls{x}{\equiv_I}$, the
$\equiv_I$-equivalence class of $x$, for any $x\in L$.

In case $L$ is sectionally complemented, this can be easily expressed
in terms of the relations of perspectivity and subperspectivity in
$L$. The following result is proved in
\cite[Theorem~III.13.20]{Birk79}.

\begin{proposition}\label{P:StIdSCL}
Let $I$ be an ideal of a sectionally complemented modular lattice~$L$.
Then~$I$ is neutral if{f} $x\sim y$ and $y\in I$ implies that
$x\in I$, for all $x$, $y\in L$.
\end{proposition}

\begin{corollary}\label{C:CharacCenCML}
Let $L$ be a sectionally complemented modular lattice. An element
$u\in L$ is neutral if{f} $x\lesssim u$ and $x\wedge u=0$ implies
that $x=0$, for all $x\in L$.
\end{corollary}

For a positive integer $n$, a \emph{homogeneous sequence of order
$n$} in a lattice~$L$ with zero is an independent (see
\cite[Definition~IV.1.9]{Grat98}) sequence $\Seq{a_0,\dots,a_{n-1}}$
of pairwise perspective elements of $L$.

The \emph{center} of a bounded lattice $L$, denoted by $\cen L$, is
the set of all complemented neutral elements of $L$. The elements of
$\cen L$ correspond exactly to the direct decompositions of
$L$. This can be expressed conveniently in the following way (see
\cite[Theorem~III.4.1]{Grat98}).

\begin{lemma}\label{L:charactCen(L)}
Let $L$ be a bounded lattice and let $a$, $b\in L$. Then the
following are equivalent:
\begin{enumerate}
\item There are bounded lattices $A$ and $B$ and an isomorphism
$f\colon L\to A\times B$ such that $f(a)=\Seq{1,0}$ and
$f(b)=\Seq{0,1}$.

\item The pair $\seq{a,b}$ is complementary in
$\cen L$, that is, $a$, $b\in\cen L$ and\linebreak
$a\oplus b=1$.
\end{enumerate}
Furthermore, $\cen(L\res a)=(\cen L)\res a$, for any $a\in\cen L$.
\end{lemma}

For the following result we refer the reader to
\cite[Theorem~III.2.9]{Grat98}.

\begin{proposition}\label{P:CenLBA}
The center of a bounded lattice $L$ is a Boolean sublattice of~$L$.
\end{proposition}

\section{Regular rings}\label{S:BasicRng}
All our rings will be associative. Most of the time they will also be
unital, with a few exceptions.
A ring $R$ is (von~Neumann) \emph{regular}, if every element $a$ of
$R$ has a \emph{quasi-inverse}, that is, an element $b$ of $R$ such
that $aba=a$. For a regular ring $R$, the set $\bL(R)$ of principal
right ideals of $R$, that is,
 \[
 \bL(R)=\setm{xR}{x\in R}=\setm{xR}{x\in R,\ x^2=x}
 \]
partially ordered by inclusion, is a sectionally complemented
modular lattice (see Section~\ref{S:BasicLatt}), with least
element $\set{0_R}$.

Hence every coordinatizable lattice is sectionally complemented and
modular. It is observed in B. J\'onsson \cite[Section~9]{Jons62} that
a bounded lattice is coordinatizable if{f} it can be coordinatized by
a regular, unital ring.

We shall need the following classical result (see
K.\,R. Goodearl \cite[Theorem~1.7]{Good91}, or
K.\,D. Fryer and I. Halperin \cite[Section~3.6]{FrHa56} for the
general, non-unital case).

\begin{proposition}\label{P:MatReg}
For any regular ring $R$ and any positive integer $n$, the ring
$\Mat_n(R)$ of all $n\times n$ matrices with entries in $R$ is
regular.
\end{proposition}

We shall need a more precise form of the result stating that $\bL(R)$
is a lattice, proved in \cite[Section~3.2]{FrHa56}.

\begin{proposition}\label{P:+capinL(R)}
Let $R$ be a regular ring and let $a$, $b\in R$ with $a^2=a$.
Furthermore, let $u$ be a quasi-inverse of $b-ab$. Then the following
statements hold:
\begin{enumerate}
\item Put $c=(b-ab)u$. Then $aR+bR=(a+c)R$.

\item Suppose that $b^2=b$ and put $d=u(b-ab)$. Then
$aR\cap bR=(b-bd)R$.
\end{enumerate}
\end{proposition}

A ring $R$ endowed with its canonical structure of right $R$-module
will be denoted by $R_R$.

\begin{corollary}\label{C:L(R)SML}
Let $R$ be a regular ring. Then $\bL(R)$ is a sectionally
complemented sublattice of $\Sub(R_R)$.
\end{corollary}

Remember that $\Sub(R_R)$ is an Arguesian lattice; hence so is
$\bL(R)$.

We shall also use the following easy consequence of
Proposition~\ref{P:+capinL(R)}, already observed in F. Micol's
thesis \cite{Mico}.

\begin{corollary}\label{C:Lfunctor}
Let $R$ and $S$ be regular rings and let $f\colon R\to S$ be a ring
homomorphism. Put $I=\ker f$. Then the following statements hold:
\begin{enumerate}
\item There exists a unique map
$g\colon\bL(R)\to\bL(S)$ such that $g(xR)=f(x)S$ for all $x\in R$.
We shall denote this map by $\bL(f)$.

\item $\bL(f)$ is a $0$-lattice homomorphism from $\bL(R)$
to $\bL(S)$.

\item There is an isomorphism $\eps$ from $\ker\bL(f)$ onto $\bL(I)$,
defined by the rule $\eps(xR)=xI$, for all $x\in I$.

\item If $f$ is a ring embedding, then $\bL(f)$ is a lattice
embedding.

\item If $f$ is surjective, then $\bL(f)$ is surjective.
\end{enumerate}
Furthermore, the correspondence $R\mapsto\bL(R)$, $f\mapsto\bL(f)$
defines a functor from the category of regular rings and ring
homomorphisms to the category of sectionally complemented modular
lattices and $0$-lattice homomorphisms. This
functor preserves direct limits.
\end{corollary}

In particular, if we identify $\setm{xR}{x\in I}$ with $\bL(I)$ (via
the isomorphism $\eps$), then we obtain the isomorphism
$\bL(R/I)\cong\bL(R)/\bL(I)$.

The following result sums up a few easy preservation
statements.

\begin{proposition}\label{P:PresCoordRP}\hfill
\begin{enumerate}
\item Any neutral ideal of a coordinatizable lattice is
coordinatizable.

\item Any homomorphic image of a coordinatizable lattice is
coordinatizable.

\item Any reduced product of coordinatizable lattices is
coordinatizable.
\end{enumerate}
\end{proposition}

\begin{proof}
(i) Let $R$ be a regular ring and let $\bI$ be a neutral ideal of
$\bL(R)$. The subset $I=\setm{x\in R}{xR\in\bI}$ is a two-sided
ideal of $R$ (see \cite[Theorem~4.3]{UnifRef}), thus, in particular,
it is a regular ring in its own right (see \cite[Lemma~1.3]{Good91}).
Furthermore, as seen above, the rule $xI\mapsto xR$
defines an isomorphism from $\bL(I)$ onto $\bI$.

(ii) In the context of (i) above,
$\bL(R)/{\bI}=\bL(R)/\bL(I)\cong\bL(R/I)$. The isomorphism
$\bL(R/I)\to\bL(R)/{\bI}$ is given by
$(\lambda+I)(R/I)\mapsto\cls{{\lambda}R}{\bI}$.

(iii) Let $\seqm{L_i}{i\in I}$ be a family of coordinatizable
lattices and let $\cF$ be a filter on~$I$. For $i\in I$, let $R_i$ be
a regular ring such that $\bL(R_i)\cong L_i$. Put
$L=\prod_{\cF}\famm{L_i}{i\in I}$ and $R=\prod_{\cF}\famm{R_i}{i\in
I}$. It is easy to verify that $\bL(R)$ is isomorphic to~$L$.
\end{proof}

The treatment of direct decompositions of a
unital ring parallels the theory for bounded lattices. For a unital
ring~$R$, we denote by
$\cen R$ the set of all central idempotents of $L$. It is well-known
that $\cen R$ is a Boolean algebra, with $a\vee b=a+b-ab$,
$a\wedge b=ab$, and $\neg a=1-a$, for all $a$, $b\in \cen R$.
The elements of $\cen R$ correspond exactly to the direct
decompositions of $R$, in a way that parallels closely
Lemma~\ref{L:charactCen(L)}.

There is a natural correspondence between the
center of a regular ring $R$ and the center of the lattice $\bL(R)$,
see \cite[Satz~VI.1.8]{FMae58}.

\begin{proposition}\label{P:CenRL}
Let $R$ be a unital regular ring. The map $e\mapsto eR$ defines an
isomorphism from $B=\cen R$ onto $\cen\bL(R)$. Furthermore,
 \[
 \cen(eR)=\setm{xR}{x\in B\res e}.
 \]
\end{proposition}

\section{Modules}\label{S:BasicModules}
A right module $E$ over a ring $R$ is \emph{semisimple}, if
the lattice $\Sub E$ of all submodules of $E$ is complemented.

\begin{proposition}\label{P:SemisCML}
Let $E$ be a semisimple right module over a unital ring $R$. Let
$S=\End E$ be the endomorphism ring of $E$, and put 
 \[
 \bI(X)=\setm{f\in R}{\im f\subseteq X},\text{ for all }
 X\in\Sub E.
 \]
Then $S$ is a regular ring and $X\mapsto\bI(X)$ defines a lattice
isomorphism from $\Sub E$ onto $\bL(S)$, with inverse the map
$fS\mapsto\im f$.
\end{proposition}

\begin{proof}
Let $f\in S$. Since $E$ is semisimple, there are submodules $X$ and
$Y$ of $E$ such that $E=X\oplus\ker f=Y\oplus\im f$. Let
$p\colon E\onto\im f$ be the projection along $Y$. For any $y\in E$,
the element $p(y)$ belongs to $\im f$, thus $p(y)=f(x)$ for a
unique element $x\in X$, that we denote by $g(y)$. Then $g\in S$ and
$f\circ g\circ f=f$, whence $S$ is regular.

Let $X\in\Sub E$. It is clear that $\bI(X)$ is a right ideal of $S$.
Furthermore, since $X$ has a direct summand in $E$, there exists a
projection $p$ of $E$ such that $\im p=X$. So,
to conclude the proof, it suffices to prove that $\bI(X)=fS$, for any
$f\in S$ with $\im f=X$. It is clear that $fS$ is contained in
$\bI(X)$. Conversely, let $g\in\bI(X)$. The submodule $\ker f$ of $E$
has a direct summand~$Y$. For any $x\in E$, the element $g(x)$
belongs to $X=\im f$, thus there exists a unique $y=h(x)$ in $Y$ such
that $g(x)=f(y)$. Hence $h\in S$ and $g=f\circ h$ belongs to $fS$;
whence $\bI(X)=fS$.
\end{proof}

In particular, we get the well-known result that for any right vector
space $E$ over any division ring, the endomorphism ring $R=\End E$
is regular and $\bL(R)\cong\Sub E$.

A nontrivial right module $E$ over a ring $R$ is \emph{simple}, if
$\Sub E=\set{\set{0},E}$. We state the classical \emph{Schur's Lemma}.

\begin{proposition}\label{P:End(simple)}
Let $E$ be a simple right module over a ring. Then $\End E$ is a
division ring.
\end{proposition}

Let a right module $E$ over a ring $R$ be expressed as a
finite direct sum $E=E_1\oplus\cdots\oplus E_n$. Let $p_i$ (resp.,
$e_i$) denote the canonical projection on $E_i$ (resp., the inclusion
map $E_i\into E$), for all $i\in\set{1,\dots,n}$.
Any endomorphism $f$ of $E$ gives raise to a system of homomorphisms
$f_{i,j}\colon E_j\to E_i$, for $i$, $j\in\set{1,\dots,n}$, defined
as $f_{i,j}=p_i\circ f\circ e_j$. Then the map
 \begin{equation}\label{Eq:BlockMatrix}
 f\mapsto\begin{pmatrix}
 f_{1,1} & \ldots & f_{1,n}\\
 \vdots & & \vdots\\
 f_{n,1} & \ldots & f_{n,n}
 \end{pmatrix}
 \end{equation}
defines an isomorphism from $\End E$ to the ring of all matrices as
in the right hand side of \eqref{Eq:BlockMatrix}, where
$f_{i,j}\in\Hom(E_j,E_i)$ for all $i$, $j\in\set{1,\dots,n}$, endowed
with canonical addition and multiplication. We shall be especially
interested in the case where all the $E_i$-s are isomorphic
submodules.

\begin{proposition}\label{P:MatrDecIsom}
In the context above, let $\gamma_i\colon E_1\to E_i$ be an
isomorphism, for all $i\in\set{1,\dots,n}$. Then the rule
 \[
 f\mapsto
\begin{pmatrix}
 \gamma_1^{-1}f_{1,1}\gamma_1 & \ldots &
 \gamma_1^{-1}f_{1,n}\gamma_n\\
 \vdots & & \vdots\\
 \gamma_n^{-1}f_{n,1}\gamma_1 & \ldots &
 \gamma_n^{-1}f_{n,n}\gamma_n
 \end{pmatrix}
 \]
defines an isomorphism from $\End E$ onto $\Mat_n(\End E_1)$.
\end{proposition}

\section{Coordinatization of lattices of length two}\label{S:CoordMn}
We denote by $M_X$ the lattice of length two and distinct atoms
$q_x$, for $x\in X$, for any nonempty set $X$. The lattices $M_3$ and
$M_7$ are diagrammed on Figure~1, Page~\pageref{Fig:2distr}. Hence
the simple lattices of length $2$ are exactly the lattices
$M_\kappa$, where $\kappa\ges3$ is a cardinal number. A
bounded lattice~$L$ has a \emph{spanning $M_X$}, if there exists a
$0,1$-lattice homomorphism $f\colon M_X\to L$
(observe that either~$f$ is one-to-one or~$L$ is trivial). The
following result is folklore.

\begin{proposition}\label{P:CoorMn}
Let $n\ges 3$ be a natural number. Then the following are equivalent:
\begin{enumerate}
\item $n-1$ is a prime power;

\item there exists a field $F$ such that $\bL(\Mat_2(F))\cong M_n$;

\item there exists a regular ring $R$ such that $\bL(R)\cong M_n$;

\item there exists a unital ring $R$ such that $\Sub(R_R)\cong M_n$;

\item there are a ring $R$ and a right $R$-module $E$ such
that $\Sub E\cong M_n$.
\end{enumerate}
\end{proposition}

\begin{proof}
(i)$\Rightarrow$(ii) Set $q=n-1$. Since the right
$\FF_q$-module $E=\FF_q\times\FF_q$ is semisimple, it follows from
Proposition~\ref{P:SemisCML} that
$\bL(\Mat_2(\FF_q))\cong\Sub E$, whence $\bL(\Mat_2(\FF_q))$ is
isomorphic to $M_n$.

(ii)$\Rightarrow$(iii) It follows from Proposition~\ref{P:MatReg}
that $R=\Mat_2(F)$ is a regular ring.

(iii)$\Rightarrow$(iv) and (iv)$\Rightarrow$(v) are trivial.

(v)$\Rightarrow$(i) By assumption,
$\Sub E=\set{\set{0},E}\cup\set{E_1,\dots,E_n}$, where
$E=E_i\oplus E_j$ for all distinct $i$, $j\in\set{1,\dots,n}$. In
particular, $E_1\oplus E_3=E_2\oplus E_3=E$, thus $E_1\cong E_2$.
Since $\Sub E\cong M_n$, the module $E$ is semisimple, whence, by
Proposition~\ref{P:SemisCML}, $S=\End E$ is a regular ring and
$\Sub E\cong\bL(S)$. Furthermore, since $E=E_1\oplus E_2$ and
$E_1\cong E_2$, it follows from Proposition~\ref{P:MatrDecIsom} that
$S\cong\Mat_2(D)$ where we put $D=\End E_1$. {}From
$\Sub E_1=\set{\set{0},E_1}$ and Proposition~\ref{P:End(simple)} it
follows that $D$ is a division ring, and hence, by using again
Proposition~\ref{P:SemisCML}, 
 \[
 M_n\cong\Sub E\cong\bL(S)\cong\bL(\Mat_2(D))\cong\Sub(D\times D).
 \]
Therefore, $D$ is a finite division ring, so the order $q$ of $D$ is a
prime power, and $n=1+q$.
\end{proof}

In particular, the first non-coordinatizable lattice of length two is
$M_7$, see Figure~1, Page~\pageref{Fig:2distr}.

By keeping track of the isomorphisms in the direction
(v)$\Rightarrow$(i) of the proof of Proposition~\ref{P:CoorMn}, we
obtain the following additional information.

\begin{proposition}\label{P:UniqCoordMn}
Let $\kappa$ be a cardinal number greater than or equal to $2$. Then
the regular rings coordinatizing $M_{1+\kappa}$ are exactly those of
the form $\Mat_2(D)$, where $D$ is a division ring with $\kappa$
elements.
\end{proposition}

In particular, for any prime power $q$, there exists exactly one
regular ring coordinatizing $M_{1+q}$, namely $\Mat_2(\FF_q)$.

We denote by~$\CC$ the class of all coordinatizable
lattices and by $\NC$ its complement
(within, say, the class of all complemented modular lattices). The
following consequence of Proposition~\ref{P:UniqCoordMn} is observed
by B. J\'onsson in the Introduction of~\cite{Jons61}.

\begin{corollary}\label{C:NonCnot1ord}
The class $\NC$ is not first-order definable. In particular, $\CC$ is
not finitely axiomatizable.
\end{corollary}

\begin{proof}
It follows from Proposition~\ref{P:CoorMn} that $M_{4k+7}$ is not
coordinatizable, for all $k<\omega$. Let $\cU$ be a
nonprincipal ultrafilter on~$\omega$. The ultraproduct, with
respect to~$\cU$, of the sequence $\Seqm{M_{4k+7}}{k<\omega}$ is
isomorphic to $M_X$, for some infinite set~$X$; thus, by
Proposition~\ref{P:UniqCoordMn}, it is coordinatizable. In
particular, the class $\NC$ is not closed under
ultraproducts, hence it is not first-order definable.
\end{proof}

\section{A first example about unions of coordinatizable lattices}
\label{S:DirUnion}

It is well-known that the center $Z(R)$ of a regular ring $R$
is regular (see \cite[Theorem~1.14]{Good91}). In particular, for each
prime~$p$, there are $a_p$, $c_p\in Z(R)$ with
 \begin{align}
 p^2a_p&=p\cdot 1_R,\label{Eq:p^2a_p=p}\\
 c_p&=1_R-pa_p.\label{Eq:cpfromap}
 \end{align}
Observe that $c_p$ is independent of the element
$a_p$ satisfying \eqref{Eq:p^2a_p=p}.

\begin{lemma}\label{L:cpinR}
The element $c_p$ is a central idempotent of $R$, for each prime
$p$. In addition, $c_pc_q=0$ for all distinct primes $p$ and $q$.
\end{lemma}

\begin{proof}
It is trivial that $c'_p=pa_p$ is idempotent; thus, so is $c_p$. As
$a_p$ is central, so are $c'_p$ and $c_p$.

Now let $p$ and $q$ be distinct primes, and put $e=c_pc_q$. {}From
$pc'_pe=p^2a_pe=pe$, it follows that
$pc_pe=0$. Since $c_pe=e$, we obtain that
$pe=0$. Similarly, $qe=0$. Since $p$ and $q$ are coprime, it
follows that $e=0$, which establishes our claim.
\end{proof}

This makes it possible to solve negatively an open problem raised by
B. J\'onsson in \cite[Section~10]{Jons62}.

\begin{proposition}\label{P:CnotclDirUn}
There exists a countable $2$-\distr\ complemented modular lattice
$L$, with a spanning $M_3$, which satisfies the two following
properties:
\begin{enumerate}
\item $L$ is a directed union of finite coordinatizable lattices;

\item $L$ is not coordinatizable.
\end{enumerate}
Consequently, the class $\CC$ of coordinatizable lattices is not
closed under countable directed unions.
\end{proposition}

\begin{proof}
Define $L$ as the set of all almost constant sequences
$x=\Seqm{x_n}{n<\omega}$ of elements of $M_4$ such that
$x(\infty)\in M_3$, endowed with componentwise
ordering. It is easy to verify that $L$ is a
countable $2$-\distr\ complemented modular lattice with a spanning
$M_3$.

For each $n<\omega$, put $L_n=(M_4)^n\times M_3$, and denote by
$f_n\colon L_n\to L$ the map defined by the rule
 \[
 f_n(\seq{x_0,\dots,x_{n-1},x})=\seq{x_0,\dots,x_{n-1},x,x,\cdots}.
 \]
Then $f_n$ is a lattice embedding from $L_n$ into $L$, and $L$ is
the increasing union of all images of the maps $f_n$. Observe that
each $L_n$ (thus each $f_n[L_n]$) is coordinatizable, see
Proposition~\ref{P:CoorMn}.

Now we prove that $L$ is not coordinatizable. Suppose, to the
contrary, that there are a regular ring $R$ and an isomorphism
$\eps\colon\bL(R)\onto L$. For all $n<\omega$, denote by
$\pi_n\colon L\onto M_4$, $x\mapsto x(n)$ the $n$-th projection, and
put $\pi_{\omega}\colon L\onto M_3$, $x\mapsto x(\infty)$.
Furthermore, put $I_n=\pi_n^{-1}\set{0}$, for all $n\les\omega$. So
$I_n$ is a neutral ideal of $L$, and, as~$L$ is a complemented
modular lattice, $\pi_n$ induces an isomorphism from $L/I_n$ onto
$\im\pi_n$. The subset $J_n=\setm{x\in R}{\eps(xR)\in I_n}$ is a
two-sided ideal of $R$, and, by Proposition~\ref{P:PresCoordRP},
we can define an isomorphism $\eps_n\colon\bL(R/J_n)\onto L/I_n$ by
the rule
 \[
 \eps_n\bigl((\lambda+J_n)(R/J_n)\bigr)=[\eps(\lambda R)]_{I_n},
 \qquad\text{for all }\lambda\in R.
 \]
In particular, for all $n<\omega$, $\bL(R/J_n)\cong M_4$, thus, by
Proposition~\ref{P:UniqCoordMn}, $R/J_n\cong\Mat_2(\FF_3)$.
Similarly, $R/J_{\omega}\cong\Mat_2(\FF_2)$.

Now we consider the central elements $a_p$, $c_p$ introduced in
\eqref{Eq:p^2a_p=p} and \eqref{Eq:cpfromap}. Projecting the equality
$4a_2=2\cdot 1_R$ on $R/J_n$, for $n\les\omega$, yields
 \begin{align}
 c_2&\in J_n,\text{ for all }n<\omega,\label{Eq:c2inInAn}\\
 c_2&\in 1+J_{\omega}.\label{Eq:c2inIom}
 \end{align}
{}From $\bigcap_{n<\omega}I_n=\set{0}$ it
follows easily that $\bigcap_{n<\omega}J_n=\set{0}$, so
\eqref{Eq:c2inInAn} yields that $c_2=0$, which contradicts
\eqref{Eq:c2inIom}.
\end{proof}

\section{Determining sequences and atomic Boolean algebras}
\label{S:DetSeq}

For models $A$ and $B$ of a first-order language $\cL$, let
$A\equiv B$ denote elementary equivalence of $A$ and $B$. A Boolean
algebra $B$ is \emph{atomic}, if the unit element of $B$ is the join
of the set $\At B$ of all atoms of $B$. The following lemma is an
immediate application of A. Tarski's classification of the complete
extensions of the theory of Boolean algebras (see C.\,C. Chang and
H.\,J. Keisler \cite[Section~5.5]{ChKe}).

\begin{lemma}\label{L:AtBaEquiv}
Let $A$ and $B$ be atomic Boolean algebras. Then $A\equiv B$
if{f}\linebreak
$\min\set{|\At A|,\aleph_0}=\min\set{|\At B|,\aleph_0}$.
\end{lemma}

Now we recall a few notions about Boolean products.
Let $\cL$ be a first-order language, let $X$ be a Boolean space,
and let $A$ be a subdirect product of a family
$\seqm{A_\fp}{\fp\in X}$ of models of $\cL$. For a first-order
formula $\varphi(\vx_0,\dots,\vx_{n-1})$ of $\cL$ and elements $a_0$,
\dots, $a_{n-1}\in A$, we put
 \[
 \bv{\varphi(a_0,\dots,a_{n-1})}=
 \setm{\fp\in X}
 {A_\fp\models\varphi(a_0(\fp),\dots,a_{n-1}(\fp))}.
 \]
We say that the subdirect product
 \begin{equation}\label{Eq:BoolProd}
 A\into\prod_{\fp\in X}A_\fp
 \end{equation}
is a \emph{Boolean product} (see \cite[Section~IV.8]{BuSa}), if the
following conditions hold:
\begin{enumerate}
\item $\bv{\varphi}$ belongs to $\Clop X$, for every \emph{atomic}
sentence $\varphi$ with parameters from~$A$;

\item for any elements $a$, $b\in A$ and any clopen subset $Y$ of
$X$, the element $a\res_Y\cup b\res_{X\setminus Y}$
belongs to $A$.
\end{enumerate}

If, in addition, the Boolean value $\bv{\varphi}$ belongs to
$\Clop X$, for \emph{every} $\cL$-sentence~$\varphi$ with
parameters from $A$, we say that \eqref{Eq:BoolProd} is a
\emph{strong Boolean product}. It is observed in M. Weese
\cite[Section~8]{Wees89} that the statement that \eqref{Eq:BoolProd}
is a strong Boolean product follows from the so-called
\emph{maximality property}, that is, for every $\cL$-formula
$\varphi(\vx,\vy_0,\dots,\vy_{n-1})$ and all $b_0$, \dots,
$b_{n-1}\in A$, there exists $a\in A$ such that
 \[
 \bv{\varphi(a,b_0,\dots,b_{n-1})}=
 \bv{\exists\vx\,\varphi(\vx,b_0,\dots,b_{n-1})}.
 \]
As the following easy lemma shows, the two notions are, in fact,
equivalent.

\begin{lemma}\label{L:EqMaxBP}
Any strong Boolean product has the maximality property.
\end{lemma}

\begin{proof}
Suppose that \eqref{Eq:BoolProd} is a strong Boolean product, let
$\varphi(\vx,\vy_0,\dots,\vy_{n-1})$ be a $\cL$-formula, and let
$b_0$, \dots, $b_{n-1}\in A$. It follows from the assumption that
$U=\bv{\exists\vx\,\varphi(\vx,b_0,\dots,b_{n-1})}$ is a clopen
subset of $X$. By definition of the $\bv{{}_{-}}$ symbol, the
equality
 \[
 U=\bigcup\famm{\bv{\varphi(x,b_0,\dots,b_{n-1})}}{x\in A}
 \]
holds, thus, since $U$ is compact, there are $k<\omega$ and elements
$a_0$, \dots, $a_{k-1}\in A$ such that
 \[
 U=\bigcup\famm{\bv{\varphi(a_j,b_0,\dots,b_{n-1})}}{j<k}.
 \]
There are pairwise disjoint clopen subsets
$U_j\subseteq\bv{\varphi(a_j,b_0,\dots,b_{n-1})}$, for $j<k$, such
that $U=\bigcup\famm{U_j}{j<k}$. Since
\eqref{Eq:BoolProd} is a Boolean product, there exists $a\in A$ such
that $a\res_{U_j}=a_j\res_{U_j}$, for all $j<k$. Therefore,
 \begin{equation}
 U=\bv{\varphi(a,b_0,\dots,b_{n-1})}.\tag*{\qed}
 \end{equation}
\renewcommand{\qed}{}
\end{proof}

The following definition is the natural extension of S. Feferman and
R.\,L. Vaught's determining sequences (see \cite[Section~6.3]{ChKe})
to strong Boolean products.

\begin{definition}\label{D:DetFam}
For a formula $\varphi$ of a first-order
language $\cL$, a pair $\seq{\Phi,\seqm{\varphi_i}{i\in I}}$ is a
\emph{determining sequence} of $\varphi$, if the following conditions
hold:
\begin{enumerate}
\item the set $I$ is finite, $\Phi$ is a first-order formula of the
language $\seq{\vee,\wedge}$ with set of free variables indexed by
$I$, and all $\varphi_i$-s are $\cL$-formulas with the same free
variables as $\varphi$;

\item $\Phi$ is \emph{isotone}, that is, the theory of Boolean
algebras infers the following statement:
 \[
 \left[\Phi\famm{\vx_i}{i\in I}\text{ and }
 \conj_{i\in I}(\vx_i\leq\vy_i)\right]
 \Longrightarrow\Phi\famm{\vy_i}{i\in I}.
 \]

\item for every strong Boolean product as in \eqref{Eq:BoolProd} and
for every $\cL$-formula $\varphi(\vec a)$ with parameters from
$A$, the following equivalence holds:
 \[
 A\models\varphi(\vec a)\Longleftrightarrow
 \Clop X\models
 \Phi\famm{\bv{\varphi_i(\vec a)}}{i\in I}).
 \]
\end{enumerate}
\end{definition}

An immediate consequence of Lemma~\ref{L:EqMaxBP} and
\cite[Theorem~8.1]{Wees89} is the following.

\begin{lemma}\label{L:ExDetSeq}
For every first-order language $\cL$, every formula of $\cL$ has a
determining sequence.
\end{lemma}

We shall use later the following application to Boolean algebras.

\begin{lemma}\label{L:ElEmbBA}
Let $A$ be a subalgebra of a Boolean algebra $B$. We suppose that
both~$A$ and $B$ are atomic, with $\At A=\At B$. Then $A$ is an
elementary submodel of~$B$.
\end{lemma}

\begin{proof}
Let $\varphi(\vx_0,\dots,\vx_{n-1})$ be a formula of the language
$\seq{\vee,\wedge}$ and let $a_0$, \dots, $a_{n-1}\in A$ such that
$A\models\varphi(\vec a)$ (where $\vec a=\seq{a_0,\dots,a_{n-1}}$);
we shall prove that $B\models\varphi(\vec a)$. Denote by~$U$ the
(finite) set of atoms of the Boolean subalgebra of $A$ generated by
$\setm{a_i}{i<n}$. We use the canonical isomorphisms
 \[
 A\cong\prod\famm{A\res u}{u\in U},\qquad
 B\cong\prod\famm{B\res u}{u\in U}.
 \]
Let $\cL$ denote the first-order language obtained by enriching the
language of Boolean algebras by $n$ additional constants
$\underline{a}_0$, \dots, $\underline{a}_{n-1}$. Let $\ol{\varphi}$
denote the sentence
$\varphi(\underline{a}_0,\dots,\underline{a}_{n-1})$ of $\cL$. The
assumption that $A\models\varphi(\vec a)$ can be rewritten as
 \begin{equation}\label{Eq:Aresumodolf}
 \prod\famm{\seq{A\res u,\vec a\wedge u}}{u\in U}\models\ol{\varphi}.
 \end{equation}
Let $u\in U$. Since $a_i\wedge u\in\set{0,u}$ for all $i<n$, it
follows from Lemma~\ref{L:AtBaEquiv} that
$\seq{A\res u,\vec a\wedge u}\equiv\seq{B\res u,\vec a\wedge u}$.
Therefore, since elementary equivalence is
preserved under direct products (see \cite[Theorem~6.3.4]{ChKe}), it
follows from \eqref{Eq:Aresumodolf} that 
 \[
 \prod\famm{\seq{B\res u,\vec a\wedge u}}{u\in U}\models\ol{\varphi},
 \]
that is, $B\models\varphi(\vec a)$.
\end{proof}

\section{Coordinatizability is not first-order}\label{S:Cnot1ord}

We put $P_p=\setm{p^{k!}}{k\in\NN}$
(where $k!=k(k-1)\cdots 2\cdot 1$), for any prime~$p$, and we put
$P=P_2\cup P_3$. We shall construct a pair of lattices $K$ and $L$.
The construction can also be performed in a similar fashion for any
pair of distinct primes, we just pick $2$ and
$3$ for simplicity. Our lattices are the following:
 \begin{align*}
 K&=\Setm{x\in\prod\famm{M_{1+k}}{k\in P}}
 {x\text{ is almost constant}};\\
 L&=\Setm{x\in\prod\famm{M_{1+k}}{k\in P}}
 {\text{both }x\res_{P_2}\text{ and }x\res_{P_3}
 \text{ are almost constant}}.
 \end{align*}
Of course, both $K$ and $L$ are $2$-\distr\ complemented
modular lattices with spanning $M_3$, and $K$ is a $0,1$-sublattice
of $L$. Furthermore, $\cen K=K\cap\two^P$ and $\cen L=L\cap\two^P$,
where $\two=\set{0,1}$. Let $\infty$, $\infty_2$, and~$\infty_3$
denote distinct objects not in~$P$. We put
 \[
 U=P\cup\set{\infty},\qquad V=P\cup\set{\infty_2,\infty_3}.
 \]
Endow $U$ with the least topology making every singleton
in $P$ clopen, and $V$ with the least topology making every
singleton of $P$ clopen as well as $P_2\cup\set{\infty_2}$ (and thus
also $P_3\cup\set{\infty_3}$). Observe that $U$ is isomorphic to
$\omega+1$ endowed with its interval topology, while $V$ is
isomorphic to the disjoint union of two copies of $U$. In particular,
both $U$ and $V$ are Boolean spaces. The \emph{canonical
map} from~$V$ onto $U$ is the map $e\colon V\to U$, whose restriction
to $P$ is the identity, and that sends both $\infty_2$ and
$\infty_3$ to $\infty$. The inverse map
$\eps\colon\Clop U\into\Clop V$,
$X\mapsto e^{-1}[X]$ is the \emph{canonical embedding} from $\Clop U$
into $\Clop V$. As an immediate application of Lemma~\ref{L:ElEmbBA},
we observe the following.

\begin{lemma}\label{L:epsEltary}
The map $\eps$ is an elementary embedding from $\Clop U$
into $\Clop V$.
\end{lemma}

Now we shall represent both $K$ and $L$ as Boolean products. We put
 \begin{align*}
 K'&=\setm{x\in\Cont(U,M_\omega)}
 {(\forall k\in P)\ x(k)\in M_{1+k}};\\
 L'&=\setm{x\in\Cont(V,M_\omega)}
 {(\forall k\in P)\ x(k)\in M_{1+k}}.
 \end{align*}
The verification of the following lemma is trivial.

\begin{lemma}\label{L:KK',LL'}
Both maps from $K'$ to $K$ and from $L'$ to $L$ defined by
restriction to~$P$ are lattice isomorphisms.
\end{lemma}

We set $M_{1+k}=M_{\omega}$ for $k\in\set{\infty,\infty_2,\infty_3}$.
With each of the lattices $K'$ and $L'$ is associated a subdirect
product, namely,
 \begin{align}
 K'&\into\prod\famm{M_{1+k}}{k\in U},
 \quad x\mapsto\seqm{x(k)}{k\in U}\label{Eq:SubdProdK}\\
 L'&\into\prod\famm{M_{1+k}}{k\in V},
 \quad x\mapsto\seqm{x(k)}{k\in V}.\label{Eq:SubdProdL}
 \end{align}
We denote by $\bv{{}_{-}}^K$ (resp., $\bv{{}_{-}}^L$) the Boolean
value function defined by the subdirect decomposition
\eqref{Eq:SubdProdK} (resp., \eqref{Eq:SubdProdL}).
We denote by $a'$ (resp., $a''$) the image of $a$ under the canonical
isomorphism from $K$ onto $K'$ (resp., from $L$ onto $L'$), for any
$a\in K$ (resp., $a\in L$).

\begin{lemma}\label{L:K'L'SBP}
Both subdirect products \textup{\eqref{Eq:SubdProdK}} and
\textup{\eqref{Eq:SubdProdL}} are strong Boolean products.
Furthermore, $\bv{\varphi(\vec a'')}^L=\eps(\bv{\varphi(\vec a')}^K)$,
for every formula $\varphi(\vx_0,\dots,\vx_{n-1})$ of
$\seq{\vee,\wedge}$ and all $a_0$, \dots, $a_{n-1}\in K$.
\end{lemma}

\begin{proof}
Let $\varphi(\vx_0,\dots,\vx_{n-1})$ be a formula of the language
$\seq{\vee,\wedge}$. An easy application of the Compactness Theorem
of first-order predicate logic gives that for any $a_0$, \dots,
$a_{n-1}\in M_\omega$, the following statements are equivalent:
\begin{itemize}
\item $M_\omega\models\varphi(a_0,\dots,a_{n-1})$;

\item
$M_{1+k}\models\varphi(a_0,\dots,a_{n-1})$ for all but
finitely many $k\in P$;

\item $M_{1+k}\models\varphi(a_0,\dots,a_{n-1})$ for
infinitely many $k\in P$.
\end{itemize}

Hence, for any finite sequence $\vec a=\seqm{a_i}{i<n}$ in $K^n$,
both Boolean values $\bv{\varphi(\vec a')}^K$ and
$\bv{\varphi(\vec a'')}^L$ are clopen, respectively in $U$ and in
$V$, and they are determined by their restrictions to $P$.
Furthermore, $\bv{\varphi(\vec a'')}^L=\eps(\bv{\varphi(\vec a')}^K)$.
\end{proof}

\begin{proposition}\label{P:KEltSubmL}
The lattice $K$ is an elementary submodel of $L$.
\end{proposition}

\begin{proof}
Let $\varphi(\vx_0,\dots,\vx_{m-1})$ be a formula of
$\seq{\vee,\wedge}$ and let $\vec a=\seqm{a_i}{i<m}\in K^m$ such that
$K\models\varphi(\vec a)$. By Lemma~\ref{L:ExDetSeq},
$\varphi$ has a determining sequence, say,
$\seq{\Phi,\seqm{\varphi_j}{j<n}}$. Since
$K'\models\varphi(\vec a')$ and by Lemma~\ref{L:K'L'SBP},
the following relation holds:
 \[
 \Clop U\models
 \Phi(\bv{\varphi_0(\vec a')}^K,\dots,
 \bv{\varphi_{n-1}(\vec a')}^K).
 \]
Hence, by Lemma~\ref{L:epsEltary},
 \[
 \Clop V\models
 \Phi(\eps(\bv{\varphi_0(\vec a')}^K),\dots,
 \eps(\bv{\varphi_{n-1}(\vec a')}^K)).
 \]
By Lemma~\ref{L:K'L'SBP},
$\eps(\bv{\varphi_j(\vec a')}^K)=\bv{\varphi_j(\vec a'')}^L$, for all
$j<n$, and hence, again by Lemma~\ref{L:K'L'SBP},
$L'\models\varphi(\vec a'')$,
and therefore $L\models\varphi(\vec a)$.
\end{proof}

\begin{proposition}\label{P:KinNC}
The lattice $K$ is not coordinatizable.
\end{proposition}

\begin{proof}
Suppose otherwise, and let $\eps\colon\bL(R)\onto K$ be an
isomorphism, where $R$ is a regular ring. We denote by
$\pi_q\colon K\onto M_{1+q}$ the canonical projection, for all
$q\in P$. The subset $I_q=\pi_q^{-1}\set{0}$ is a neutral ideal of
$K$, and, as $K$ is a complemented modular lattice, $\pi_q$ induces
an isomorphism from $K/I_q$ onto $M_{1+q}$. The subset
$J_q=\setm{x\in R}{\eps(xR)\in I_q}$ is a two-sided ideal of $R$,
and, by Proposition~\ref{P:PresCoordRP}, we can define an
isomorphism $\eps_q\colon\bL(R/J_q)\onto K/I_q$ by the rule
 \[
 \eps_q\bigl((\lambda+J_q)(R/J_q)\bigr)=[\eps(\lambda R)]_{I_q},
 \qquad\text{for all }\lambda\in R.
 \]
In particular, for all $q\in P$, $\bL(R/J_q)\cong M_{1+q}$, thus, by
Proposition~\ref{P:UniqCoordMn},
 \begin{equation}\label{Eq:R/JqcongM2Fq}
 R/J_q\cong\Mat_2(\FF_q).
 \end{equation}
Now we consider again the
elements $a_p$ and $c_p$ introduced in \eqref{Eq:p^2a_p=p} and
\eqref{Eq:cpfromap}. It follows from Proposition~\ref{P:CenRL} that
$u_p=\eps(c_pR)$ belongs to the center of $K$. {}From $c_2c_3=0$
(see Lemma~\ref{L:cpinR}) it follows that $u_2\wedge u_3=0$. As
$\cen K$ consists of all almost constant elements of $\set{0,1}^P$,
it follows that either $u_2(\infty)=0$ or
$u_3(\infty)=0$. Suppose, for example, that $u_2(\infty)=0$. In
particular, there exists $q\in P_2$ such that $\pi_q(u_2)=0$.
As $q$ is a power of $2$ and by \eqref{Eq:R/JqcongM2Fq}, we get
$2a_2\in J_q$, so
 \begin{equation}\label{Eq:c2in1+Iq}
 c_2\in 1+J_q.
 \end{equation}
On the other hand,
 \begin{align*}
 0&=[u_2]_{I_q}&&(\text{because }\pi_q(u_2)=0)\\
 &=[\eps(c_2R)]_{I_q}&&(\text{by the definition of }u_2)\\
 &=\eps_q\bigl((c_2+J_q)(R/J_q)\bigr)
 &&(\text{by the definition of }\eps_q),
 \end{align*}
thus, as $\eps_q$ is an isomorphism, $c_2\in J_q$, which contradicts
\eqref{Eq:c2in1+Iq}.
\end{proof}

\begin{proposition}\label{P:LinC}
The lattice $L$ is coordinatizable.
\end{proposition}

\begin{proof}
It is obvious that $L\cong L_2\times L_3$, where we put
 \[
 L_p=\Setm{x\in\prod\famm{M_{1+q}}{q\in P_p}}
 {x\text{ is almost constant}},\qquad\text{for each prime }p.
 \]
Hence it suffices to prove that $L_p$ is coordinatizable, for each
prime $p$.

Put $S_q=\Mat_2(\FF_q)$, for each prime power $q$.
As $\FF_{p^{k!}}$ is a subfield of $\FF_{p^{(k+1)!}}$ for each
$k<\omega$, we can define a unital ring $R_p$ by
 \[
 R_p=\Setm{\lambda\in\prod\famm{S_q}{q\in P_p}}{\lambda
 \text{ is almost constant}}.
 \]
It is easy to verify that $R_p$ is a regular ring.
We shall prove that $L_p\cong\bL(R_p)$. Fix a one-to-one enumeration
$\seqm{\xi_k}{k\in\NN}$ of $\FF_{p^{\infty}}$ such that
$\FF_{p^{n!}}=\setm{\xi_k}{1\les k\les p^{n!}}$ for all $n<\omega$.
We put $\alpha_0=\begin{pmatrix}0&0\\ 0&1\end{pmatrix}$ and
$\alpha_k=\begin{pmatrix}1&0\\ \xi_k&0\end{pmatrix}$, for each
$k\in\NN$. For all $q\in P_p$, there exists a
unique isomorphism $\eta_q\colon\bL(S_q)\onto M_{1+q}$ such that
 \[
 \eta_q(\alpha_kS_q)=q_k\quad(\text{the }k
 \text{-th atom of }M_{1+q}),
 \quad\text{for all }k\in\set{0,1,\dots,q}.
 \]
We can define a map
$\eps_p\colon\bL(R_p)\to\prod\famm{M_{1+q}}{q\in P_p}$ by the rule
 \[
 \eps_p(\lambda R_p)=\seqm{\eta_q(\lambda_qS_q)}{q\in P_p},
 \quad\text{for all }\lambda\in R_p.
 \]
For any $\lambda\in R_p$, there exists $m\in P_p$ such that
$\lambda_q=\lambda_m$ for all $q\ges m$ in $P_p$. If $\lambda_m$ is
neither zero nor invertible in $S_m$, then there exists a unique
$k\in\set{0,1,\dots,m}$ such that $\lambda_mS_m=\alpha_kS_m$, thus
$\lambda_qS_q=\alpha_kS_q$ for all $q\ges m$ in $P_p$, and thus
$\eps_p(\lambda R_p)$ is almost constant (with limit $q_k$). This
holds trivially in case $\lambda_m$ is either zero or invertible,
therefore the range of $\eps_p$ is contained in $L_p$.
Now it follows from Corollary~\ref{C:Lfunctor}(ii) that
$\eps_p$ is a lattice homomorphism from~$\bL(R_p)$ onto~$L_p$.

For idempotent $\alpha$, $\beta\in R_p$, if
$\eps_p(\alpha R_p)=\eps_p(\beta R_p)$, then
$\eta_q(\alpha_qS_q)=\eta_q(\beta_qS_q)$ for all $q\in P_p$, thus
(as the $\eta_q$s are isomorphisms) $\alpha_q=\beta_q\alpha_q$ and
$\beta_q=\alpha_q\beta_q$ for all $q\in P_p$, so
$\alpha=\beta\alpha$ and $\beta=\alpha\beta$, and so
$\alpha R_p=\beta R_p$. Therefore, $\eps_p$ is one-to-one.

Let $x\in L_p$. If $x(\infty)\in\set{0,1}$, then, as each $\eta_q$
is an isomorphism, there exists $\lambda\in R_p$, with limit either
$0$ or $1$, such that $\eps_p(\lambda R_p)=x$. Now suppose that
$x(\infty)=q_k$, with $k<\omega$. There exists $m\ges k$ in $P_p$
such that $x_q=q_k$ holds for all $q\ges m$ in $P_p$. For each $q<m$
in $P_p$, there exists $\lambda_q\in S_q$ such that
$\eta_q(\lambda_qS_q)=x_q$. Put $\lambda_q=\alpha_k$, for all
$q\ges m$ in $P_p$. Then $\lambda\in R_p$ and
$\eps_p(\lambda R_p)=x$. Therefore, the map~$\eps_p$ is surjective,
and so it is an isomorphism.
\end{proof}

By combining Propositions~\ref{P:KEltSubmL}, \ref{P:KinNC}, and
\ref{P:LinC}, we obtain a negative solution to J\'onsson's Problem.

\begin{theorem}\label{T:CNCNot1}
Neither the class $\CC$ of coordinatizable lattices nor its
complement~$\NC$ are first-order classes. In fact, there are
countable, $2$-\distr\ lattices $K$ and $L$ with spanning $M_3$ such
that $K$ is an elementary sublattice of $L$, the lattice~$K$ belongs
to $\NC$, and the lattice $L$ belongs to~$\CC$.
\end{theorem}

By using sheaf-theoretical methods, we could prove that \emph{every
countable $2$-distributive complemented modular lattice with a
spanning $M_{\omega}$ is coordinatizable}. Hence the use of prime
numbers in Theorem~\ref{T:CNCNot1} is somehow unavoidable. As we
shall see in the next section, this result does not extend to the
uncountable case.

\section{An uncountable non-coordinatizable lattice with a spanning
$M_{\omega}$}\label{S:UnctbleNonCoord}

We start with an elementary lemma of linear algebra.

\begin{lemma}\label{L:EmbProjGeom}
Let $E$ be a unital ring, let $F$ be a division ring, let $n\in\NN$,
and let $\varphi\colon\Mat_n(E)\to\Mat_n(F)$ be a unital ring
homomorphism. There are a unital ring homomorphism
$\sigma\colon E\to F$ and a matrix $a\in\GL_n(F)$ such that
 \[
 \varphi(x)=a(\sigma x)a^{-1},\text{ for all }x\in\Mat_n(E),
 \]
where $\sigma x$ denotes the matrix obtained by applying $\sigma$ to
all the entries of $x$.
\end{lemma}

\begin{proof}
Put $R=\Mat_n(E)$ and $S=\Mat_n(F)$.
Let $\seqm{e_{i,j}^E}{1\les i,j\les n}$ denote the canonical system
of matrix units of $R$, and similarly for $S$. Then
$\seqm{\varphi(e_{i,j}^E)}{1\les i,j\les n}$ is a system of matrix
units of $S$, thus, since $F$ is a division ring, there
exists $a\in\GL_n(F)$ such that $\varphi(e_{i,j}^E)=ae_{i,j}^Fa^{-1}$
for all $i$, $j\in\set{1,\dots,n}$. Hence, conjugating $\varphi$ by
$a^{-1}$, we reduce the problem to the case where
$\varphi(e_{i,j}^E)=e_{i,j}^F$, for all $i$, $j\in\set{1,\dots,n}$.
Fix $i$, $j\in\set{1,\dots,n}$ and $x\in E$. As
$e_{i,j}^Ex=e_{i,i}^Exe_{j,j}^E$, the value
$\varphi(e_{i,j}^Ex)$ belongs to $e_{i,i}^FSe_{j,j}^F$, thus it has
the form $e_{i,j}^F\sigma_{i,j}(x)$, for a unique
$\sigma_{i,j}(x)\in F$.

It is obvious that all $\sigma_{i,j}$-s are unit-preserving additive
homomorphisms from $E$ to~$F$. Applying the ring homomorphism
$\varphi$ to the equalities
 \[
 e_{i,i}^Ex=e_{i,j}^E(e_{j,i}^Ex)=(e_{i,j}^Ex)e_{j,i}^E,
 \text{ for all }x\in E,
 \]
we obtain that $\sigma_{i,i}=\sigma_{i,j}=\sigma_{j,i}$. Thus
$\sigma_{i,j}=\sigma_{1,1}$, for all $i$, $j\in\set{1,\dots,n}$.
Denote this map by $\sigma$. So $\varphi(x)=\sigma x$, for all $x\in
R$, and $\sigma$ is a ring homomorphism.
\end{proof}

\begin{corollary}\label{C:EmbProjGeom}
Let $E$ and $F$ be division rings and let
$\varphi\colon\Mat_2(E)\into\Mat_2(F)$ be a unital ring embedding. If
$\varphi$ is not an isomorphism, then the complement of the image
of $\bL(\varphi)$ in $\bL(\Mat_2(F))$ has cardinality at least $|E|$.
\end{corollary}

\begin{proof}
By using Lemma~\ref{L:EmbProjGeom}, we may conjugate $\varphi$ by a
suitable element of $\Mat_2(F)$ to reduce the problem to the case
where $\varphi(x)=\sigma x$ identically on $\Mat_2(E)$, for a
suitable unital ring embedding $\sigma\colon E\into F$. Since
$\varphi$ is not an isomorphism, $\sigma$ is not surjective. In
particular, $|F\setminus\sigma E|\ges|E|$. Now observe that the
matrices of the form
$e_{\lambda}=\begin{pmatrix}1 & 0\\ \lambda & 0 \end{pmatrix}$, with
$\lambda\in F$, are idempotent matrices of $\Mat_2(F)$
with pairwise distinct image spaces. Furthermore,
$e_{\lambda}\Mat_2(F)$ does not belong to the image of
$\bL(\varphi)$, for all $\lambda\in F\setminus\sigma E$.
\end{proof}

Now our counterexample. For an infinite cardinal number $\kappa$,
whose successor cardinal we denote by $\kappa^+$, we put
 \begin{align}
 \ol{L}_{\kappa}&=
 \Setm{x\in(M_{\kappa+1})^{\kappa^+}}{x\text{ is almost constant}},
 \label{Eq:DefolLkap}\\
 L_{\kappa}&=\Setm{x\in\ol{L}_{\kappa}}{x(\infty)\in M_{\kappa}},
 \label{Eq:DefLkap}
 \end{align}
both ordered componentwise. It is obvious that $L_{\kappa}$
is a $2$-\distr\ complemented modular lattice with a
spanning~$M_{\kappa}$.

\begin{theorem}\label{T:LnotCoord}
The lattice $L_{\kappa}$ is not coordinatizable, for any infinite
cardinal number $\kappa$.
\end{theorem}

\begin{proof}
Otherwise let $R$ be a unital regular ring and let
$\eps\colon\bL(R)\onto L_{\kappa}$ be an isomorphism. For
$x\in M_{\kappa+1}$ and $i<\kappa^+$, let $x\cdot u_i$ denote the
element of $L_{\kappa}$ with
$i$-th component $x$ and all other components zero; furthermore, put
$u_i=1\cdot u_i$. There are central idempotents $a_i\in\cen R$ such
that $\eps(a_iR)=u_i$, for all $i<\kappa^+$. Since
$\bL(a_iR)\cong L_{\kappa}\res u_i\cong M_{\kappa}$
(see Proposition~\ref{P:PresCoordRP}(i)), there exists, by
Proposition~\ref{P:UniqCoordMn}, a division ring
$E_i$ of cardinality $\kappa$ such that $a_iR\cong\Mat_2(E_i)$ (as
rings). Put
 \[
 J_X=\bigoplus\famm{a_iR}{i\in X},
 \text{ for all }X\subseteq\kappa^+,
 \]
and put $J=J_{\kappa^+}$. Observe that $I=\eps\bL(J)$ is the ideal of
all almost null elements of~$L_{\kappa}$. It follows that
 \begin{equation}\label{Eq:R/Jtowiso}
 \bL(R/J)\cong\bL(R)/\bL(J)\cong L_{\kappa}/I\cong M_{\kappa},
 \end{equation}
with the canonical isomorphism
$\zeta\colon\bL(R/J)\onto L_{\kappa}/I$ of
Proposition~\ref{P:PresCoordRP}(ii) given by
 \[
 \zeta((x+J)R/J)=\cls{\eps(xR)}{I},\text{ for all }x\in R.
 \]
Furthermore, it follows from \eqref{Eq:R/Jtowiso} and
Proposition~\ref{P:UniqCoordMn} that $R/J\cong\Mat_2(E)$,
for some division ring $E$ with $\kappa$ elements.
In particular, $R/J$ has $\kappa$ elements. For any $\lambda\in R/J$,
pick $\dot{\lambda}\in R$ such that $\lambda=\dot{\lambda}+J$. Of
course, we may take $\dot{0}_{R/J}=0_R$ and
$\dot{1}_{R/J}=1_R$. For $\alpha$, $\beta$, $\gamma\in R/J$ such that
$\gamma=\alpha-\beta$, there exists a finite sub\-set~$X$ of
$\kappa^+$ such that
$\dot{\gamma}\equiv\dot{\alpha}-\dot{\beta}\pmod{J_X}$. By doing the
same for the product map $\seq{\alpha,\beta}\mapsto\alpha\beta$,
the zero, and the unit of $R/J$ and forming the union of all
corresponding $X$-s, we obtain a subset $X$ of $\kappa^+$ of
cardinality at most $\kappa$ such that
$p_i\colon\lambda\mapsto a_i\dot{\lambda}$ defines a unital ring
homomorphism from $R/J$ to $a_iR$, for all $i\in\kappa^+\setminus X$.
Since $R/J$ is simple, $p_i$ is an embedding.

Put $\dot{x}_{\lambda}=\eps(\dot{\lambda}R)$ (an element of
$L_{\kappa}$), for all $\lambda\in R/J$. Observe that
$\zeta(\lambda(R/J))=\cls{\dot{x}_{\lambda}}{I}$, for all
$\lambda\in R/J$; in particular
 \begin{equation}\label{Eq:DescrFullL/I}
 L_{\kappa}/I=\setm{\cls{\dot{x}_{\lambda}}{I}}{\lambda\in R/J}.
 \end{equation}
For $\alpha$, $\beta$, $\gamma\in R/J$ such that
$\cls{\dot{x}_{\gamma}}{I}=
\cls{\dot{x}_{\alpha}}{I}\vee\cls{\dot{x}_{\beta}}{I}$,
there exists a finite subset $Y$ of~$\kappa^+$ such that the equality
$\dot{x}_{\gamma}(i)=\dot{x}_{\alpha}(i)\vee\dot{x}_{\beta}(i)$
holds for all $i\in\kappa^+\setminus Y$. By doing the same for the
meet and the constants $0$ and $1$, and then taking
the union of~$X$ and all corresponding $Y$-s, we obtain a
subset~$Y$ of~$\kappa^+$ containing $X$, with at most $\kappa$
elements, such that
$g_i\colon\cls{\dot{x}_{\lambda}}{I}\mapsto
\dot{x}_{\lambda}\wedge u_i$ defines a
$\seq{\vee,\wedge,0,1}$-homomorphism from
$L_{\kappa}/I$ into $L_{\kappa}\res u_i$, for all
$i\in\kappa^+\setminus Y$; since
$L_{\kappa}/I$ is simple, $g_i$ is, actually, an embedding.

Let $\zeta_i\colon\bL(a_iR)\to L_{\kappa}\res u_i$,
$xR\mapsto\eps(xR)$ denote the canonical isomorphism, for all
$i<\kappa^+$. We shall verify that the following diagram is
commutative, for $i\in\kappa^+\setminus Y$.
 \[
\def\labelstyle{\displaystyle}
\xymatrix{
\bL(R/J)\ar[rr]^{\bL(p_i)}\ar[d]_{\zeta} &&
\bL(a_iR)\ar[d]^{\zeta_i}\\
L_{\kappa}/I\ar[rr]_{g_i} && L_{\kappa}\res u_i
}
 \]
Let $\lambda\in R/J$. It is immediate that
$\zeta_i\circ\bL(p_i)(\lambda(R/J))=\eps(a_i\dot{\lambda}R)$. On the
other hand, we compute
 \begin{align*}
 g_i\circ\zeta(\lambda(R/J))&=g_i(\cls{\dot{x}_{\lambda}}{I})\\
 &=\dot{x}_{\lambda}\wedge u_i\\
 &=\eps(\dot{\lambda}R)\wedge\eps(a_iR)\\
 &=\eps(\dot{\lambda}R\wedge a_iR)\\
 &=\eps(a_i\dot{\lambda}R)&&(\text{because }a_i\in\cen R),
 \end{align*}
which completes the verification of the commutativity of the diagram
above. By applying \eqref{Eq:DescrFullL/I} to the classes modulo $I$
of constant functions, we obtain that for all $q\in M_{\kappa}$,
there exists $\lambda_q\in R/J$ such that the set
$Z_q=\setm{i\in\kappa^+}{\dot{x}_{\lambda_q}(i)\neq q}$
is finite; whence
 \begin{equation}\label{Eq:reachallMo}
 g_i(\cls{\dot{x}_{\lambda_q}}{I})=q\cdot u_i,\text{ for all }
 i\in\kappa^+\setminus(Y\cup Z_q).
 \end{equation}
Furthermore, there exists a subset $Z$ of $\kappa^+$
containing $Y\cup\bigcup\famm{Z_q}{q\in M_{\kappa}}$,
with at most $\kappa$ elements, such that
$\dot{x}_{\lambda}$ is constant on $\kappa^+\setminus Z$, with
value, say, $y_{\lambda}\in M_{\kappa}$, for all $\lambda\in R/J$.
Hence,
 \begin{equation}\label{Eq:Avoidqo}
 g_i(\cls{\dot{x}_{\lambda}}{I})=\dot{x}_{\lambda}\wedge u_i
 =y_{\lambda}\cdot u_i\neq q_{\kappa}\cdot u_i,
 \text{ for all }\seq{\lambda,i}\in(R/J)\times(\kappa^+\setminus Z).
 \end{equation}
Therefore, by \eqref{Eq:reachallMo} and \eqref{Eq:Avoidqo},
we obtain that
 \[
 \im g_i=(L_{\kappa}\res u_i)\setminus\set{q_{\kappa}\cdot u_i},
 \text{ for all }i\in\kappa^+\setminus Z.
 \]
Since both maps~$\zeta$ and $\zeta_i$ are isomorphisms and the
diagram above is commutative, the complement in $\bL(a_iR)$ of the
range of $\bL(p_i)$ is also a singleton, for all
$i\in\nobreak\kappa^+\setminus\nobreak Z$. Since $R/J\cong\Mat_2(E)$
and $a_iR\cong\Mat_2(E_i)$, we obtain, by
Corollary~\ref{C:EmbProjGeom}, a contradiction.
\end{proof}

Pushing the argument slightly further yields the following strong
negative statement.

\begin{theorem}\label{T:NotLinf}
There is no formula $\theta$ of $\cL_{\infty,\infty}$ such that the
class of $2$-distributive coordinatizable lattices is the class
of all models of~$\theta$.
\end{theorem}

\begin{proof}
For any division ring~$D$ with infinite cardinal $\kappa$, the ring
of all almost constant $\kappa^+$-sequences of elements of
$\Mat_2(D)$ coordinatizes the lattice $\ol{L}_{\kappa}$ defined
in~\eqref{Eq:DefolLkap}; whence
$\ol{L}_{\kappa}\in\CC$. We have seen in Theorem~\ref{T:LnotCoord}
that $L_{\kappa}\in\NC$. Of course,~$L_{\kappa}$ is a sublattice of
$\ol{L}_{\kappa}$. Since $\kappa$ is arbitrarily large, it is
sufficient, in order to conclude the proof, to establish
that~$L_{\kappa}$ is a $\cL_{\kappa,\kappa}$-elementary submodel of
$\ol{L}_{\kappa}$.

So we need to prove that $\ol{L}_{\kappa}\models\varphi$ implies that
$L_{\kappa}\models\varphi$, for every
$\seq{\vee,\wedge}$-sentence~$\varphi$ in $\cL_{\kappa,\kappa}$ with
parameters from $L_{\kappa}$. The only nontrivial instance of the
proof is to verify that
$\ol{L}_{\kappa}\models\exists\vec\vx\,\psi(\vec a,\vec\vx)$ implies
that $L_{\kappa}\models\exists\vec\vx\,\psi(\vec a,\vec\vx)$, for
every formula $\psi$ in $\cL_{\kappa,\kappa}$ for which we have
already proved elementariness, with a list of parameters
$\vec a=\seqm{a_{\xi}}{\xi<\alpha}$ from $L_{\kappa}$ and a list of
free variables $\vec\vx=\seqm{\vx_{\eta}}{\eta<\beta}$, where
$\alpha$, $\beta<\kappa$. So let us fix a list
$\vec b=\seqm{b_{\eta}}{\eta<\beta}$ from $\ol{L}_{\kappa}$ such
that $\ol{L}_{\kappa}\models \psi(\vec a,\vec b)$. Since
$\alpha$, $\beta<\kappa$, there are $\gamma<\kappa^+$ and an
automorphism $\sigma$ of $M_{\kappa+1}$ such that the following
statements hold:
\begin{enumerate}
\item $x(\zeta)=x(\infty)$, for all
$x\in\setm{a_{\xi}}{\xi<\alpha}\cup\setm{b_{\eta}}{\eta<\beta}$ and
all $\zeta\in\kappa^+\setminus\gamma$;

\item $\sigma(a)=a$, for all
$a\in\setm{a_{\xi}(\infty)}{\xi<\alpha}$;

\item $\sigma(b)\in M_{\kappa}$, for all
$b\in\setm{b_{\eta}(\infty)}{\eta<\beta}$.
\end{enumerate}
Denote by $\tau$ the automorphism of $\ol{L}_{\kappa}$ defined by
the rule
 \[
 \tau(x)(\zeta)=\begin{cases}
 x(\zeta),&\text{if }\zeta<\gamma,\\
 \sigma(x(\zeta)),&\text{if }\gamma\les\zeta,
 \end{cases}
 \qquad\text{for all }\seq{x,\zeta}\in\ol{L}_{\kappa}\times\kappa^+.
 \]
Then $\tau$ fixes all $a_{\xi}$-s while the element
$c_{\eta}=\tau(b_{\eta})$ belongs to $L_{\kappa}$, for all
$\eta<\beta$. {}From $\ol{L}_{\kappa}\models\psi(\vec a,\vec b)$ it
follows that $\ol{L}_{\kappa}\models\psi(\vec a,\vec c)$, thus,
by the induction hypothesis,
$L_{\kappa}\models\psi(\vec a,\vec c)$, and therefore
$L_{\kappa}\models\exists\vec\vx\,\psi(\vec a,\vec\vx)$.
\end{proof}

\section{Appendix: Large partial three-frames are finitely
axiomatizable}\label{S:LP3Frme}

For a positive integer $n$ and a bounded lattice $L$, we say that $L$
has a \emph{large partial $n$-frame}, if there exists a homogeneous
sequence $\Seq{a_0,\dots,a_{n-1}}$ of order $n$ in $L$ such that~$L$ is
generated by $a_0$ as a neutral ideal. It is clear that the existence
of a large partial $(n+1)$-frame implies the existence of a large
partial $n$-frame.

Having a large partial $3$-frame does not appear to be
a first-order condition \emph{a priori}. However, we shall now prove
that it is.

\begin{proposition}\label{P:LP3FRM1}
Let $L$ be a complemented Arguesian lattice. Then $L$ has a large
partial $3$-frame if{f} there are $a_0$, $a_1$, $a_2$, $b\in L$ such
that
\begin{enumerate}
\item $a_0\oplus a_1\oplus a_2\oplus b=1$;

\item $a_i\sim a_j$, for all distinct $i$, $j<3$;

\item $b\lesssim a_0\oplus a_1$.
\end{enumerate}
In particular, for a complemented Arguesian lattice, having a large
partial $3$-frame can be expressed by a single first-order sentence.
\end{proposition}

\begin{proof}
It is obvious that the given condition implies that
$\Seq{a_0,a_1,a_2}$ is a homogeneous sequence such that the neutral
ideal generated by $a_0$ is $L$.

Conversely, suppose that $L$ has a large partial $3$-frame.
We shall make use of the \emph{dimension monoid} $\Dim L$ of
$L$ introduced in \cite{WDim}. As in \cite{WDim}, we
denote by $\DD(x,y)$ the element of $\Dim L$ representing the abstract
``distance'' between elements~$x$ and~$y$ of
$L$. Since $L$ has a zero, we put $\DD(x)=\DD(0,x)$, for all $x\in L$.
We shall also use the result, proved in
\cite[Theorem~5.4]{WDim}, that the dimension monoid of a modular
lattice is a refinement monoid.

Putting $\eps=\DD(1)$ and applying the unary function $\DD$ to the
parameters of a large partial $3$-frame of~$L$, we obtain that there
are $n\in\NN$ and $\alpha$, $\beta\in\Dim L$ such that the
following relations hold:
 \begin{align}
 3\alpha+\beta&=\eps;\label{Eq:3alphbeteps}\\
 \beta&\leq n\alpha.\label{Eq:betnalph}
 \end{align}
Furthermore, by J\'onsson's Theorem, $L$ is coordinatizable, thus
\emph{normal} as defined in~\cite{WDim}. This implies easily the
following statement:
 \begin{equation}\label{Eq:DelNorm}
 (\DD(x)=\DD(y)\text{ and }x\wedge y=0)\Longrightarrow
 x\sim y,\text{ for all }x,\,y\in L.
 \end{equation}
Since $\Dim L$ is a refinement monoid,
\eqref{Eq:betnalph} implies (see \cite[Lemma~3.1]{WDim}) the existence
of elements $\alpha_k\in\Dim L$, for $0\les k\les n$, such that
 \[
 \alpha=\sum_{0\les k\les n}\alpha_k\qquad\text{and}\qquad
 \beta=\sum_{0\les k\les n}k\alpha_k.
 \]
Put $\ol{\alpha}=\sum_{0\les k\les n}
\bigl(\intp{\frac{k}{3}}+1\bigr)\alpha_k$ and
$\ol{\beta}=\sum_{0\les k\les n}
\bigl(k-3\intp{\frac{k}{3}}\bigr)\alpha_k$, where $\intp{x}$ denotes
the largest integer below $x$, for every rational number $x$. Hence
we immediately obtain
 \begin{equation}\label{Eq:3alb+bb}
 3\ol{\alpha}+\ol{\beta}=\eps.
 \end{equation}
To prove that $\ol{\beta}\leq2\ol{\alpha}$, it suffices to prove that
$k-3\intp{\frac{k}{3}}\les2\intp{\frac{k}{3}}+2$ for all
$k\in\set{0,1,\dots,n}$, which is immediate.

Since the $\DD$ function is a V-measure (see
\cite[Corollary~9.6]{WDim}), there are $a_0$, $a_1$, $a_2$, $b\in L$
such that $a_0\oplus a_1\oplus a_2\oplus b=1$ in $L$ while
$\DD(a_i)=\ol{\alpha}$ for all $i<3$ and $\DD(b)=\ol{\beta}$.
Hence, by \eqref{Eq:DelNorm}, $\Seq{a_0,a_1,a_2}$ is a homogeneous
sequence. Furthermore,
$\DD(b)=\ol{\beta}\leq2\ol{\alpha}=\DD(a_0\oplus a_1)$, thus, by
\cite[Corollary~9.4]{WDim}, there are $x_0\leq a_0$, $x_1\leq a_1$,
and $b_0$, $b_1\leq b$ such that $\DD(x_0)=\DD(b_0)$,
$\DD(x_1)=\DD(b_1)$, and $b=b_0\oplus b_1$. It follows again
from \eqref{Eq:DelNorm} that $b=b_0\oplus b_1\sim x_0\oplus x_1$,
whence $b\lesssim a_0\oplus a_1$.
\end{proof}

We remind the reader of J\'onsson's
Extended Coordinatization Theorem (cf.
Page~\pageref{T:JonssThm}), which states that for complemented
Arguesian lattices, existence of a large partial $3$-frame implies
coordinatizability. In particular, lattices with a large partial
$3$-frame are not enough to settle J\'onsson's Problem.

\section{Open problems}\label{S:OpenPbs}

Some of our problems will be formulated in the language of
descriptive set theory. We endow the powerset $\Pow(X)\cong\two^X$
with the product topology of the discrete topological space
$\two=\set{0,1}$, for any set~$X$. So $\Pow(X)$ is compact
Hausdorff, metrizable in case $X$ is countable. Hence the space
$\BS=\Pow(\omega^2)\times\Pow(\omega^3)\times\Pow(\omega^3)$, endowed
with the product topology, is also compact metrizable. We endow it
with its canonical recursive presentation (see Y.\,N.
Moschovakis \cite{Mosc}).

We define $\BL$ as the set of all triples $\xi=\seq{E,M,J}\in\BS$
such that $E$ is a partial ordering on a nonzero initial
segment $m$ of $\omega$ on which $J$ and $M$ are,
respectively, the join and the meet operation with respect to $E$,
and the lattice $L_\xi=\seq{m,E,M,J}$ is complemented modular.

Since stating that a structure is a
complemented modular lattice can be expressed by a
finite set of $\forall\exists$ axioms, it is not hard to verify that
$\BL$ is a $\Pi^0_2$ subset of $\BS$. Put
 \[
 \CL=\setm{\xi\in\BL}{L_\xi\text{ is coordinatizable}},
 \]
the set of \emph{real codes} of coordinatizable lattices.
As $\CL$ is defined by a second-order existential statement, it is a
$\Sigma_1^1$ subset of $\BS$.

\begin{problem}\label{Pb:CoordSig11}
Is $\CL$ a Borel subset of $\BS$?
\end{problem}

\begin{problem}\label{Pb:Orthoc}
Is the set of real codes of countable complemented modular lattices
admitting an orthocomplementation a Borel subset of $\BS$?
\end{problem}

By using sheaf-theoretical methods, we could prove that the analogue
of Problem~\ref{Pb:Orthoc} for $2$-distributive lattices with
spanning $M_3$ has a \emph{positive} solution. In fact, the obtained
condition is first-order.

\begin{problem}\label{Pb:UpwEltary}
Let $K$ be an elementary sublattice of a countable bounded lattice
$L$. If $K$ is coordinatizable, is $L$ coordinatizable?
\end{problem} 

Our next problem is related to a possible weakening of the
definition of coordinatizability.

\begin{problem}\label{Pb:WeakCoord}
If the join-semilattice $\Subc M$ of all finitely generated
submodules of a right module $M$ is a complemented lattice, is it
coordinatizable?
\end{problem}

Proposition~\ref{P:SemisCML} answers Problem~\ref{Pb:WeakCoord}
positively only in case $M$ is n{\oe}therian.

\begin{problem}\label{Pb:Dec2CML}
Describe the elementary invariants of $2$-\distr\ complemented
modular lattices. Is the theory of $2$-\distr\ complemented modular
lattices decidable?
\end{problem}

\begin{problem}\label{Pb:Fin2Fin}
If a finite lattice $L$ can be embedded into some complemented
modular lattice, can it be embedded into some \emph{finite}
complemented modular lattice?
\end{problem}

A complemented modular lattice $L$ is \emph{uniquely
coordinatizable}, if there exists a unique (up to isomorphism)
regular ring $R$ such that $L\cong\bL(R)$.

\begin{problem}\label{Pb:UniqCoord}
Is the class of uniquely coordinatizable lattices a first-order
class?
\end{problem}

It is established in \cite{Jons60} that every complemented Arguesian
lattice with a large partial $3$-frame (see Section~\ref{S:LP3Frme})
is uniquely coordinatizable. The uniqueness part is
\cite[Theorem~9.4]{Jons60}.

\section*{Acknowledgment}

Part of this work was done during the author's visit at the TU
Darmstadt in November 2002. The hospitality of the Arbeitsgruppe~14
and Christian Herrmann's so inspiring coaching on modular lattices
are highly appreciated.

\end{document}